\documentclass[11pt]{article}
\usepackage{amscd}
\usepackage{amsfonts}
\usepackage{amsmath}
\usepackage{amssymb}
\usepackage{amsthm}
\usepackage{bbm}
\usepackage{CJK}
\usepackage{fancyhdr}
\usepackage{graphicx}
\usepackage{hyperref}
\usepackage{indentfirst}
\usepackage{latexsym}
\usepackage{mathrsfs}
\usepackage{xypic}

\newtheorem{theorem}{Theorem}[section]
\newtheorem{lemma}[theorem]{Lemma}
\newtheorem{definition}[theorem]{Definition}
\newtheorem{proposition}[theorem]{Proposition}
\newtheorem{example}[theorem]{Example}

\newtheorem{corollary}[theorem]{Corollary}

\usepackage[top=1in,bottom=1in,left=1.25in,right=1.25in]{geometry}
\def\<{\langle}
\def\>{\rangle}
\def\a{\alpha}
\def\b{\beta}

\def\c{\cdot}

\def\g{\gamma}

\date{}
\begin{document}
\renewcommand{\baselinestretch}{1.2}
\renewcommand{\arraystretch}{1.0}
\title{\bf On split regular BiHom-Leibniz superalgebras}
\author{{\bf Shuangjian Guo$^{1}$,  Shengxiang Wang$^{2}$\footnote
        { Corresponding author(Shengxiang Wang):~~wangsx-math@163.com} }\\
{\small 1. School of Mathematics and Statistics, Guizhou University of Finance and Economics} \\
{\small  Guiyang  550025, P. R. of China} \\
{\small 2.~ School of Mathematics and Finance, Chuzhou University}\\
 {\small   Chuzhou 239000,  P. R. of China}}
 \maketitle
\begin{center}
\begin{minipage}{13.cm}

{\bf \begin{center} ABSTRACT \end{center}}
The goal of this paper is to study the structure of split regular BiHom-Leiniz superalgebras, which is a natural
generalization of split regular Hom-Leiniz algebras and  split regular BiHom-Lie superalgebras.  By developing techniques of connections
of roots for this kind of algebras, we show that such a split regular BiHom-Leiniz superalgebras $\mathfrak{L}$ is of the form  $\mathfrak{L}=U+\sum_{\a}I_\a$ with $U$ a subspace of a maximal abelian
subalgebra $H$ and any $I_{\a}$, a well described ideal of $\mathfrak{L}$, satisfying $[I_\a, I_\b]= 0$ if
$[\a]\neq [\b]$. In the case of  $\mathfrak{L}$ being of maximal length,  the simplicity of $\mathfrak{L}$ is also
characterized in terms of connections of roots.


 \end{minipage}
 \end{center}
 \normalsize\vskip1cm

\section*{INTRODUCTION}
\def\theequation{0. \arabic{equation}}
\setcounter{equation} {0}

The notion of Hom-Lie algebras was first introduced by Hartwig,
Larsson and Silvestrov in \cite{Hartwig}, who developed an approach to deformations 
of the Witt and  Virasoro algebras basing on $\sigma$-deformations.
In fact, Hom-Lie algebras include Lie algebras as a subclass, but the deformation  of Lie algebras twisted by a homomorphism.

Leibniz superalgebra is a sort of non-commutative generalization of the Lie superalgebra structure,
with non skew-symmetric bracket, where the Jacobi superidentity is replaced by the so called
Leibniz superidentity. The notion of Leibniz superalgebras was firstly introduced in \cite{Abdaoui17}, although graded
Leibniz algebra was considered before in work \cite{Livernet98}. As Leibniz algebras are a generalization
of Lie algebras \cite{Loday93}, then many of the features of Leibniz superalgebras are
generalization of Lie superalgebras. As a generalization of Hom-Lie superalgebras,   Hom-Leibniz superalgebras are introduced and  presented the methods to construct these superalgebras  in \cite{Wang15}. The cohomology of Hom-Leibniz superalgebras is studied in \cite{Abdaoui14}.   some characterizations of Hom-Leibniz superalgebras are given and some of their basic properties are found, the existence of a Hom-Lie-Yamaguti superalgebra structure on any (multiplicative) left Hom-Leibniz superalgebra is proved in \cite{Gaparayi18}.

A BiHom-algebra is an algebra in such a way that the identities defining the structure
are twisted by two homomorphisms $\phi$ and $\psi$. This class of algebras was introduced from a
categorical approach in \cite{Graziani} which as an extension of the class of Hom-algebras. If the two
linear maps are the same automorphisms, BiHom-algebras will be return to Hom-algebras.
These algebraic structures include BiHom-associative algebras, BiHom-Lie algebras and
BiHom-bialgebras. The representation theory of BiHom-Lie algebras was introduced by
Cheng and Qi in \cite{Cheng16}, in which, BiHom-cochain complexes, derivation, central extension, derivation
extension, trivial representation and adjoint representation of BiHom-Lie algebras were
studied. More applications of  BiHom-algebras, BiHom-Lie superalgebras,  BiHom-Lie colour algebras and BiHom-Novikov algebras
can be found in (\cite{Liu17}, \cite{Wang16}, \cite{Abdaoui17}, \cite{Guo2018}).

The class of the split algebras is specially related to addition quantum numbers, graded contractions and deformations. For instance, for a physical system
which displays a symmetry, it is interesting to know the detailed structure of the split
decomposition, since its roots can be seen as certain eigenvalues which are the additive
quantum numbers characterizing the state of such system. Determining the structure of
split algebras will become more and more meaningful in the area of research in mathematical physics. Recently, in (\cite{Albuquerque 2018}-\cite{Cao18}, \cite{Guo2019}), the structure of different classes of split algebras have been determined by the techniques of connections of roots. The purpose
of this paper is to  consider the structure of split regular BiHom-Leiniz superalgebras by the
techniques of connections of roots based on some work in \cite{Cao18} and  \cite{Zhang2018}.

This paper is organized
as follows.  In Section 2, we prove that such an arbitrary split regular BiHom-Leiniz superalgebras $\mathfrak{L}$ is of the form  $\mathfrak{L}=U+\sum_{\a}I_\a$ with $U$  a subspace of a maximal abelian subalgebra $H$ and any $I_{\a}$, a well described ideal of $\mathfrak{L}$, satisfying $[I_\a, I_\b]= 0$ if
$[\a]\neq [\b]$.  In Section 3, we present that under certain conditions, in the case of $\mathfrak{L}$ being of maximal
length, the simplicity of the algebra is characterized.

\section{Preliminaries}
\def\theequation{\arabic{section}.\arabic{equation}}
\setcounter{equation} {0}

Throughout this paper, we will denote by $\mathbb{N}$ the set of all nonnegative integers and by $\mathbb{Z}$ the set of all integers.  Split regular BiHom-Leiniz superalgebras  are considered of arbitrary dimension and over an arbitrary base field  $\mathbb{K}$. And we recall some basic definitions and results related to our paper from \cite{Cao18} and  \cite{Zhang2018}.

\noindent{\bf 1.1. Hom-Leiniz algebra }
A Hom-Leiniz algebra $\mathfrak{L}$ is an algebra $\mathfrak{L}$  over an arbitrary base filed $\mathbb{K}$, endowed with a bilinear
product
\begin{eqnarray*}
[\c,\c]: \mathfrak{L}\times \mathfrak{L}\rightarrow \mathfrak{L},
\end{eqnarray*}
and a homomorphism $\phi: \mathfrak{L}\rightarrow \mathfrak{L}$
\begin{eqnarray*}
&&[\phi(x), [y, z]]=[[x,y],\phi(z)]+[\phi(y), [x,z]]~~\mbox{(Hom-Leibniz~ identity)}
\end{eqnarray*}
holds for any $x, y, z\in \mathfrak{L}$.

If $\phi$ is furthermore an algebra automorphism, that is,  a linear bijective on such that $\phi([x,y])=[\phi(x),\phi(y)]$ for any $x,y\in \mathfrak{L}$, then $\mathfrak{L}$ is called a regular Hom-Leiniz algebra.

\noindent{\bf 1.2. BiHom-Lie superalgebra }  A BiHom-Lie superalgebra $L$ is a $\mathbb{Z}_2$-graded algebra $L=L_{\overline{0}}\oplus L_{\overline{1}}$, endowed with an even bilinear mapping $[\c,\c]: L\times L\rightarrow L$ and two homomorphisms $\phi, \psi: L\rightarrow L$
 satisfying the following conditions, for all  $x\in L_{\overline{i}},y\in L_{\overline{j}},z\in L_{\overline{k}}$ and $\overline{i},\overline{j},\overline{k}\in \mathbb{Z}_2$:
\begin{eqnarray*}
&& [x,y]\subset L_{\overline{i}+\overline{j}},\\
&&\phi\circ\psi=\psi\circ\phi,\\
&&[\psi(x),\phi(y)]=-(-1)^{\overline{i}\overline{j}}[\psi(y),\phi(x)],\\
&&(-1)^{\overline{k}\overline{i}}[\psi^{2}(x),[\psi(y),\phi(z)]]+(-1)^{\overline{i}\overline{j}}[\psi^{2}(y),[\psi(z),\phi(x)]]+(-1)^{\overline{j}\overline{k}}[\psi^{2}(z),[\psi(x),\phi(y)]]=0.
\end{eqnarray*}

When $\phi$ and $\psi$ are algebra automorphisms, it is said  that $L$ is a regular BiHom-Lie superalgebra.

\section{Decomposition}
\def\theequation{\arabic{section}. \arabic{equation}}
\setcounter{equation} {0}

\begin{definition}
A BiHom-Leiniz superalgebra $\mathfrak{L}$ is a $\mathbb{Z}_2$-graded algebra $\mathfrak{L}=\mathfrak{L}_{\bar{0}}\oplus \mathfrak{L}_{\bar{1}}$  over an arbitrary base filed $\mathbb{K}$, endowed with a bilinear
product
\begin{eqnarray*}
[\c,\c]: \mathfrak{L}\times \mathfrak{L}\rightarrow \mathfrak{L},
\end{eqnarray*}
and two superspace homomorphisms $\phi, \psi: \mathfrak{L}\rightarrow \mathfrak{L}$
\begin{eqnarray*}
&&[x, y]\subset \mathfrak{L}_{\bar{i}+\bar{j}}\\
&&\phi\circ\psi=\psi\circ\phi,\\
&&[\phi\psi(x), [y, z]]=[[\psi(x),y],\psi(z)]+(-1)^{\bar{i}\bar{j}}[\psi(y), [\phi(x),z]]~~\mbox{(BiHom-Leibniz~ superidentity)}
\end{eqnarray*}
holds for any $x\in \mathfrak{L}_{\bar{i}}, y\in \mathfrak{L}_{\bar{j}},z\in \mathfrak{L}_{\bar{k}}$ and $\bar{i},\bar{j},\bar{k}\in \mathbb{Z}_2$.
\end{definition}
If $\phi,  \psi$ are furthermore algebra automorphisms, it is said that  $\mathfrak{L}$ is called a regular BiHom-Leiniz superalgebra.

 \begin{example}
 Let $(\mathfrak{L}, [\cdot, \cdot])$ be a Leiniz superalgebra and $\phi, \psi: \mathfrak{L}\rightarrow \mathfrak{L}$ are two commuting automorphisms of  Leiniz  superalgebras. If we endow the underlying linear space $\mathfrak{L}$ with a new product $[\cdot,\cdot]': \mathfrak{L}\times \mathfrak{L}\rightarrow \mathfrak{L}$ defined by $[x,y]'=[\phi(x),\psi(y)]$ for any $x,y\in \mathfrak{L}$, we know that $(\mathfrak{L},  [\cdot, \cdot]',  \phi, \psi)$ becomes a regular BiHom-Leiniz superalgebra.
\end{example}

\begin{example}
Consider  the 5-dimensional $\mathbb{Z}_2$-graded vector space $\mathfrak{L}=\mathfrak{L}_{\bar{0}}\oplus \mathfrak{L}_{\bar{1}}$,  over an arbitrary base filed $\mathbb{K}$ of characteristic different from 2, with basis $\{u_1, u_2,u_3\}$ of $\mathfrak{L}_{\bar{0}}$ and $\{e_1, e_2\}$ of $\mathfrak{L}_{\bar{1}}$; and where the nonzero products on these elements are induced by the following relations:
\begin{eqnarray*}
&&[u_2, u_1]=-u_3, ~~[u_1, u_2]=u_3,~~ [u_1, u_3]=-2u_1,\\
&&[u_3, u_1]=2u_1, ~~[u_3, u_2]=-2u_2,~~ [u_2, u_3]=2u_2,\\
&&[e_1, u_2]=e_2,~~[e_1, u_3]=-e_1,~~[e_2, u_1]=e_1,~~[e_2, u_3]=e_2.
\end{eqnarray*}
Then by considering the superspace homomorphism
\begin{eqnarray*}
\phi, \psi: \mathfrak{L} \rightarrow \mathfrak{L}
\end{eqnarray*}
defined by
\begin{eqnarray*}
\phi(u_i)=u_i, ~~~i=1,2,3, ~~\phi(e_1)=-e_1,~~~\phi(e_2)=-e_2.\\
\psi(u_i)=-u_i, ~~~i=1,2,3, ~~\psi(e_1)=e_1,~~~\psi(e_2)=e_2.
\end{eqnarray*}
Then we have that $\mathfrak{L}=\mathfrak{L}_{\bar{0}}\oplus \mathfrak{L}_{\bar{1}}$  becomes a  split regular BiHom-Leiniz superalgebra.
\end{example}

Note that $\mathfrak{L}_{\bar{0}}$ is a BiHom-Leiniz algebra. Moreover, if the identity $[x, y]=-(-1)^{\bar{i}\bar{j}}[y, x]$ holds, we have that  BiHom-Leibniz superidentity becomes BiHom-Jacobi superidentity, and so BiHom-Leibniz superalgebras generalize both BiHom-Leibniz algebras and BiHom-Lie superalgebras. The usual regularity concepts will be understood in the graded sense. That is, a supersubalgebra $A=A_{\bar{0}}\oplus A_{\bar{1}}$ of $\mathfrak{L}$ is a graded  subspace $A=A_{\bar{0}}\oplus A_{\bar{1}}$  satisfying
\begin{eqnarray*}
[A,A]\subset \mathfrak{L},~~~ \phi(A)=A, ~~~\psi(A)=A.
\end{eqnarray*}
  An ideal $I$ of $\mathfrak{L}$ is a graded  subspace $I=I_{\bar{0}}\oplus I_{\bar{1}}$ of $\mathfrak{L}$ such that
\begin{eqnarray*}
[I,\mathfrak{L}]+[\mathfrak{L}, I]\subset I, ~~~~\phi(I)=I,  ~~~\psi(I)=I.
\end{eqnarray*}
The graded ideal $\mathfrak{J}$ generated by
\begin{eqnarray*}
\{[x, y]+{-1}^{\bar{i}\bar{j}}[y, x]: x\in \mathfrak{L}_{\bar{i}}, y\in \mathfrak{L}_{\bar{j}}, \bar{i},\bar{j}\in \mathbb{Z}_2\}
\end{eqnarray*}
plays an important role in the theory since it determines the non-super Lie character of $\mathfrak{L}$. From BiHom-Leibniz superidentity, it is straightforward to check that this ideal satisfies
\begin{eqnarray}
[\mathfrak{L}, \mathfrak{J}]=0.
\end{eqnarray}

Let us introduce the class of split algebras in the framework of BiHom-Leibniz superalgebras in a similar
way to the cases of BiHom-Lie superalgebras. Denote by $H=H_{\bar{0}}\oplus H_{\bar{1}}$ a maximal abelian subalgebra $H$ of a BiHom-Leiniz superalgebra $\mathfrak{L}$.  For a linear functional
\begin{eqnarray*}
\a:H_{\bar{0}}\rightarrow \mathbb{K},
\end{eqnarray*}
we define the root space of $\mathfrak{L}$ associated to $\a$ as the subspace
\begin{eqnarray*}
\mathfrak{L}_{\a}:=\{v_{\a}\in \mathfrak{L}:[h_{\bar{0}},\phi(v_{\a})]=\a(h)\phi\psi(v_{\a}), \mbox{for any $h_{\bar{0}}\in H_{\bar{0}}$}\}.
\end{eqnarray*}
The elements $\a:H_{\bar{0}}\rightarrow \mathbb{K}$ satisfying $\mathfrak{L}_{\a}\neq 0$ are called roots of $\mathfrak{L}$ with respect to $H$ and we denote $\Lambda:=\{\a\in (H_{\bar{0}})^{\ast}/\{0\}:\mathfrak{L}_{\a}\neq 0\}$. We call that $\mathfrak{L}$ is a split regular BiHom-Leibniz superalgebra with respect to $H$ if
\begin{eqnarray*}
\mathfrak{L}=H\oplus \bigoplus_{\a\in \Lambda}\mathfrak{L}_{\a}.
\end{eqnarray*}
We also say that $\Lambda$ is the root system of $\mathfrak{L}$.

\begin{example}
Let $\mathfrak{L}=H\oplus (\bigoplus_{\a\in \Gamma}\mathfrak{L}_{\a})$ be a Leibniz superalgebra, $\phi, \psi: \mathfrak{L}\rightarrow \mathfrak{L}$ two commuting  automorphisms such that $\phi(H)=H, \psi(H)=H$. By Example  2.2, we know that $(\mathfrak{L}, [\c,\c]', \phi, \psi)$ is a regular BiHom-Leibniz superalgebra. Then we have
\begin{eqnarray*}
\mathfrak{L}=H\oplus (\bigoplus_{\a\in \Gamma}\mathfrak{L}_{\a\psi^{-1}})
\end{eqnarray*}
makes of the regular BiHom-Leibniz superalgebra $(\mathfrak{L},  [\c,\c]', \phi, \psi)$ being the roots system $\Lambda=\{\a\psi^{-1}:\a\in \Gamma\}$.
\end{example}

\begin{example}
By Example  2.3, we know that $(\mathfrak{L}, [\c,\c], \phi, \psi)$ is a regular BiHom-Leibniz superalgebra. Then we have
\begin{eqnarray*}
\mathfrak{L}=H\oplus \mathfrak{L}_{\a}\oplus \mathfrak{L}_{-\a}\oplus \mathfrak{L}_{\b}\oplus \mathfrak{L}_{-\b}
\end{eqnarray*}
where $H=<u_3>,  \mathfrak{L}_{\a}=<u_2>, \mathfrak{L}_{-\a}=<u_1>, \mathfrak{L}_{\b}=<e_1>$ and    $\mathfrak{L}_{-\b}=<e_2>$, being $\a,\b: H\rightarrow \mathbb{K}$ defined by
$\a(\lambda u_3)=2\lambda$ and $\b(\lambda u_3)=-\lambda, \lambda\in  \mathbb{K}$.
\end{example}
To simplify notation, the mappings $\phi|_H, \psi|_H, \phi^{-1}|_H, \psi^{-1}|_H: H\rightarrow H$ will be denoted by $\phi, \psi $ and $\phi^{-1}, \psi^{-1}$ respectively.

The following two lemmas are  analogous  to the results  of \cite{Zhang2018}.

\begin{lemma}
 Let $(\mathfrak{L}, [\cdot, \cdot], \phi, \psi)$ be a split regular BiHom-Leiniz superalgebra. Then,
for any $\a,\b\in \Lambda\cup \{0\}$,

(1) $\phi(\mathfrak{L}_{\a})=\mathfrak{L}_{\a\phi^{-1}}$ and $\phi^{-1}(\mathfrak{L}_{\a})=\mathfrak{L}_{\a\phi}$,

(2) $\psi(\mathfrak{L}_{\a})=\mathfrak{L}_{\a\psi^{-1}}$ and $\psi^{-1}(\mathfrak{L}_{\a})=\mathfrak{L}_{\a\psi}$,

(3) $[\mathfrak{L}_{\a},\mathfrak{L}_{\b}]\subset \mathfrak{L}_{\a\phi^{-1}+\b\psi^{-1}}$,
\end{lemma}
\begin{lemma}
The following assertions hold

(1) If $\a\in \Lambda$, then $\a \phi^{-z_1}\psi^{-z_2}\in \Lambda$ for any $z_1,z_2 \in \mathbb{Z}$,

(2) $\mathfrak{L}_{0}=H$.
\end{lemma}

In what follows, $\mathfrak{L}$ denotes a split regular BiHom-Leiniz superalgebra and
\begin{eqnarray*}
\mathfrak{L}=H\oplus (\bigoplus_{\a\in \Lambda}\mathfrak{L}_{\a})
\end{eqnarray*}
the corresponding root spaces decomposition. Given a linear functional $\a:H_{\bar{0}}\rightarrow \mathbb{K}$, we denote by $-\a:H_{\bar{0}}\rightarrow \mathbb{K} $ the element in $H^{\ast}_{\bar{0}}$ defined by $(-\a)(h_{\bar{0}}):=-\a(h_{\bar{0}})$. We write
\begin{eqnarray*}
-\Lambda:=\{-\a:\a\in \Lambda\} ~~~\mbox{and}~~ \pm \Lambda:\Lambda\cup (-\Lambda).
\end{eqnarray*}

\begin{definition}
 Let $\a,\b\in \Lambda$. We will say that $\a$ is connected to $\b$ if either
\begin{eqnarray*}
\b=\epsilon\a\phi^{z_1}\psi^{z_2}~~\mbox{for some $z_1,z_2\in \mathbb{Z}$ and $\epsilon\in \{-1,1\}$}
\end{eqnarray*}
or there exists $\{\a_1,\cdot\cdot\cdot,\a_k\}\subset \pm\Lambda$ with $k\geq 2$, such that \\
1. $\a_{1}\in \{\a \phi^{-n}\psi^{-r}:~~n, r\in \mathbb{N}\}$.\\
2. $\a_1\phi^{-1}+\a_2\psi^{-1}\in \pm\Lambda$,

 $\a_1\phi^{-2}+\a_2\phi^{-1}\psi^{-1}+\a_3\psi^{-1}\in \pm\Lambda$,

 $\a_1\phi^{-3}+\a_2\phi^{-2}\psi^{-1}+\a_3\phi^{-1}\psi^{-1}+\a_4\psi^{-1}\in \pm\Lambda$,

 $\cdot\cdot\cdot\cdot\cdot$

 $\a_1\phi^{-i}+\a_2\phi^{-i+1}\psi^{-1}+\a_3\phi^{-i+2}\psi^{-1}+\cdot\cdot\cdot+\a_i\phi^{-1}\psi^{-1}+\a_{i+1}\psi^{-1}\in \pm\Lambda$,

 $\a_1\phi^{-k+2}+\a_2\phi^{-k+3}\psi^{-1}+\a_3\phi^{-k+4}\psi^{-1}+\cdot\cdot\cdot+\a_{k-2}\phi^{-1}\psi^{-1}+\a_{k-1}\psi^{-1}\in \pm\Lambda$.\\
 3. $\a_1\phi^{-k+1}+\a_2\phi^{-k+2}\psi^{-1}+\a_3\phi^{-k+3}\psi^{-1}+\cdot\cdot\cdot+\a_{i}\phi^{-k+i}\psi^{-1}+\cdot\cdot\cdot+\a_{k-1}\phi^{-1}\psi^{-1}+\a_{k}\psi^{-1} \in \{\pm \b\phi^{-m}\psi^{-s}:m, s\in \mathbb{N}\}$.

 We will also say that  $\{\a_1,\cdot\cdot\cdot,\a_k\}$ is a connection from $\a$ to $\b$.
\end{definition}
The proof of the next result is analogous to the one of \cite{Cao18}.

\begin{proposition}
The relation $\sim$ in $\Lambda$, defined by $\alpha\sim \b$ if and only if  $\a$ is connected to $\b$, is an equivalence relation.
\end{proposition}

  By  Proposition 2.9 we can consider the quotient set
  \begin{eqnarray*}
 \Lambda/\sim=\{[\a]:\a\in \Lambda\},
  \end{eqnarray*}
with $[\a]$ being the set of nonzero roots which are connected to $\a$.
Our next goal is to associate an ideal $I_{[\a]}$ to $[\a]$. Fix $[\a]\in \Lambda/\sim$, we start by defining
\begin{eqnarray*}
I_{H,[\a]}=span_{\mathbb{K}}\{[\mathfrak{L}_{\b\psi^{-1}}, \mathfrak{L}_{-\b\phi^{-1}}]:\b\in [\a]\}.
\end{eqnarray*}
Now we define
\begin{eqnarray*}
V_{[\a]}:=\bigoplus_{\b\in [\a]}\mathfrak{L}_{\b}.
\end{eqnarray*}
Finally, we denote by $I_{[\a]}$ the direct sum of the two subspaces above:
\begin{eqnarray*}
I_{[\a]}:=I_{H,[\a]}\oplus V_{[\a]}.
\end{eqnarray*}

\begin{proposition} For any $[\a]\in \Lambda/\sim$, the following assertions hold.

(1). $[I_{[\a]},I_{[\a]}]\subset I_{[\a]}$,

(2). $\phi(I_{[\a]})=I_{[\a]}$ and  $\psi(I_{[\a]})=I_{[\a]}$, 

(3). For any $[\b]\neq [\a]$, we have $[I_{[\a]}, I_{[\b]}]=0$.
\end{proposition}
{\bf Proof.} (1) First we check that $[I_{[\a]},I_{[\a]}]\subset I_{[\a]}$, we can write
\begin{eqnarray}
 [I_{[\a]},I_{[\a]}]&=&[I_{H,[\a]}\oplus V_{[\a]},I_{H,[\a]}\oplus V_{[\a]}]\nonumber\\
 &\subset&[I_{H,[\a]}, V_{[\a]}]+[V_{[\a]}, I_{H,[\a]}]+[V_{[\a]}, V_{[\a]}],
\end{eqnarray}
Given $\b\in [\a]$, we have $[I_{H,[\a]}]\subset \mathfrak{L}_{\b\psi^{-1}}$. Since $\b\psi^{-1}\in [\a]$, it follows that $[I_{H,[\a]}, \mathfrak{L}_{\b}]\subset V_{[\a]}$.

By a similar argument, we get $[\mathfrak{L}_{\b}, I_{H,\a}]\subset V_{[\a]}$.

Next we consider  $[V_{[\a]},V_{[\a]}]$. If we take $\b,\g\in [\a]$ and $\bar{i}, \bar{j}\in \mathbb{Z}_2$ such that $[\mathfrak{L}_{\b, \bar{i}}, \mathfrak{L}_{\g}, \bar{j}]\neq 0$, then $[\mathfrak{L}_{\b}, \mathfrak{L}_{\g}]\subset \mathfrak{L}_{\b\phi^{-1}+\g\psi^{-1},  \bar{i}+\bar{j}}$.
If $\b\phi^{-1}+\g\psi^{-1}=0$, we get $[\mathfrak{L}_{\b}, \mathfrak{L}_{-\g}]\subset H$ and so $[\mathfrak{L}_{\b}, \mathfrak{L}_{-\g}]\subset I_{H,[\a]}$. Suppose that $\b\phi^{-1}+\g\psi^{-1}\in \Lambda$. We infer that $\{\b,\g\}$ is a connection from $\b$ to $\b\phi^{-1}+\g\psi^{-1}$. The transitivity of $\sim$ now gives that $\b\phi^{-1}+\g\psi^{-1}\in [\a]$ and so $[\mathfrak{L}_{\b}, \mathfrak{L}_{\g}]\subset V_{[\a]}$. Hence
\begin{eqnarray}
[V_{[\a]},V_{[\a]}]\in I_{[\a]}.
\end{eqnarray}
 From (2.1) and (2.2), we get $[I_{[\a]},I_{[\a]}]\subset I_{[\a]}$.

(2) It is easy to check that $\phi(I_{[\a]})=I_{[\a]}$ and $\psi(I_{[\a]})=I_{[\a]}$.

(3) We will study the expression $[I_{[\a]},I_{[\b]}]$. Notice that
\begin{eqnarray}
[I_{[\a]},I_{[\b]}]&=&[I_{H,[\a]}\oplus V_{[\a]},I_{H,[\b]}\oplus V_{[\b]}]\nonumber\\
 &\subset&[I_{H,[\a]}, V_{[\b]}]+[ V_{[\a]},I_{H,[\b]} ]+[V_{[\a]}, V_{[\b]}],
\end{eqnarray}

First we  consider $[V_{[\a]}, V_{[\b]}]$ and suppose that there exist $\a_1\in [\a], \b_1\in [\b]$ and  $\bar{j}, \bar{k}\in \mathbb{Z}_2$such that $[\mathfrak{L}_{\a_1, \bar{j}}, \mathfrak{L}_{\b_1, \bar{k}}]\neq 0$. As necessarily $\a_1\phi^{-1}\neq-\b_1\psi^{-1}$, then $\a_1\phi^{-1}+\b_1\psi^{-1}\in \Lambda$. So $\{\a_1,\b_1, -\a_1\phi^{-2}\psi\}$ is a connection between $\a_1$ and $\b_1$. By the transitivity of the connection relation we see $\a\in [\b]$, a contradicition. Hence $[\mathfrak{L}_{\a_1, \bar{j}}, \mathfrak{L}_{\b_1, \bar{k}}]=0$, and so
\begin{eqnarray}
[V_{[\a]}, V_{[\b]}]=0.
\end{eqnarray}
Next we consider the first summand $[I_{H,[\a]}, V_{[\b]}]$ on the right hand side of (2.4), and suppose that there exist $\a_1\in [\a]$ and $\b_1\in [\b]$ such that
\begin{eqnarray*}
[ \mathfrak{L}_{\b_1}, [\mathfrak{L}_{\a_1}, \mathfrak{L}_{-\a_1}]]= 0.
\end{eqnarray*}
If
\begin{eqnarray*}
[ \mathfrak{L}_{\b_1}, [\mathfrak{L}_{\a_1}, \mathfrak{L}_{-\a_1}]]\neq 0,
\end{eqnarray*}
then BiHom-Leibniz~  superidentity gives
\begin{eqnarray*}
0&\neq& [ \phi\psi(\mathfrak{L}_{\b_1}), [\mathfrak{L}_{\a_1}, \mathfrak{L}_{-\a_1}]]\\
&\subset&[[\psi(\mathfrak{L}_{\b_1}), \mathfrak{L}_{\a_1}],\psi(\mathfrak{L}_{-\a_1})]+[\psi(\mathfrak{L}_{\a_1}), [\phi(\mathfrak{L}_{\b_1}), \mathfrak{L}_{-\a_1}]].
\end{eqnarray*}
Hence
\begin{eqnarray*}
[\phi(\mathfrak{L}_{\b_1}), \mathfrak{L}_{\a_1}]+[\phi(\mathfrak{L}_{\b_1}), \mathfrak{L}_{-\a_1}]\neq 0,
\end{eqnarray*}
which contradicts (2.5). Therefore, $[\phi\psi(\mathfrak{L}_{\b_1}), [\mathfrak{L}_{\a_1}, \mathfrak{L}_{-\a_1}]]=0$.

 Consequently, $[I_{H,[\a]},V_{[\b]}]=0.$   In a similar way we can prove  $[ V_{[\a]},I_{H,[\b]} ]=0$ and the proof is complete. \hfill $\square$

\begin{definition}
A BiHom-Leiniz superalgebra $\mathfrak{L}$ is called simple if $[\mathfrak{L}, \mathfrak{L}]\neq 0$  and its only ideals are $\{0\}, \mathfrak{J}$ and $\mathfrak{L}$.
\end{definition}

\begin{theorem}
The following assertions hold

(1) For any $[\a]\in \Lambda/\sim$, the linear space $I_{[\a]}=I_{H,[\a]}+V_{[\a]}$ of $\mathfrak{L}$ associated to $[\a]$ is an ideal of $\mathfrak{L}$.

(2) If $\mathfrak{L}$ is simple, then there exists a connection from $\a$ to $\b$ for any $\a,\b\in \Lambda$ and $H=\sum_{\a\in \Lambda}([\mathfrak{L}_{\a\psi^{-1}}, \mathfrak{L}_{-\a\phi^{-1}}])$.
\end{theorem}
{\bf Proof.} (1) Since $[I_{[\a]},H]+[H, I_{[\a]}]\subset I_{[\a]}$, by  Lemmas 2.6 and  2.7, we have
\begin{eqnarray*}
[I_{[\a]}, \mathfrak{L}]=[I_{[\a]}, H\oplus(\bigoplus_{\b\in [\a]} \mathfrak{L}_{\b})\oplus (\bigoplus_{\g\notin [\a]} \mathfrak{L}_{\g})]\subset I_{[\a]}
\end{eqnarray*}
and
\begin{eqnarray*}
[ \mathfrak{L}, I_{[\a]},]=[ H\oplus(\bigoplus_{\b\in [\a]} \mathfrak{L}_{\b})\oplus (\bigoplus_{\g\notin [\a]} \mathfrak{L}_{\g}), I_{[\a]}]\subset I_{[\a]}. 
\end{eqnarray*}
According to Propositions 2.9 and  2.10,  we have
\begin{eqnarray*}
&&[I_{[\a]}, \mathfrak{L}]+[\mathfrak{L}, I_{[\a]}]\\
&=&[I_{[\a]}, (H\oplus(\bigoplus_{\b\in [\a]} \mathfrak{L}_{\b})\oplus (\bigoplus_{\g\notin [\a]} \mathfrak{L}_{\g}))]\\
&&+[(H\oplus(\bigoplus_{\b\in [\a]} \mathfrak{L}_{\b})\oplus (\bigoplus_{\g\notin [\a]} \mathfrak{L}_{\g})), I_{[\a]}]\subset I_{[\a]}.
\end{eqnarray*}
As we also have $\phi(I_{[\a]})=I_{[\a]}$ and $\psi(I_{[\a]})=I_{[\a]}$. So we conclude that $I_{[\a]}$ is an ideal of $\mathfrak{L}$.

(2) The simplicity of $\mathfrak{L}$ implies $I_{[\a]}\in \{\mathfrak{J}, \mathfrak{L}\}$. If some $\a \in \Lambda$ is such that $I_{[\a]}= \mathfrak{L}$, then $[\a]=\Lambda$ and $H=\sum_{\a\in \Lambda}([\mathfrak{L}_{\a\psi^{-1}}, \mathfrak{L}_{-\a\phi^{-1}}])$.  Otherwise, if $I_{[\a]}= \mathfrak{J}$ for any $\a \in \Lambda$, then $[\a]=[\b]$ for any $\a,\b\in \Lambda$, and so $[\a]=\Lambda$.  Thus $H=\sum_{\a\in \Lambda}([\mathfrak{L}_{\a\psi^{-1}}, \mathfrak{L}_{-\a\phi^{-1}}])$.
\hfill $\square$

\begin{theorem} We have
\begin{eqnarray*}
\mathfrak{L}=U+\sum_{[\a]\in \Lambda/\sim}I_{[\a]},
\end{eqnarray*}
where $U$ is a linear complement in $H$ of $span_{\mathbb{K}}\{[\mathfrak{L}_{\a\psi^{-1}}, \mathfrak{L}_{-\a\phi^{-1}}] :\a\in \Lambda \}$ and any $I_{[\a]}$ is one of the ideals of $\mathfrak{L}$ described in Theorem 2.12, satisfying $[I_{[\a]},I_{[\b]}]=0$ if $[\a]\neq[\b]$.
\end{theorem}
{\bf Proof.}  $I_{[\a]}$ is well defined and is an ideal of $\mathfrak{L}$, being clear that
 \begin{eqnarray*}
 \mathfrak{L}=H\oplus\sum_{[\a]\in \Lambda} \mathfrak{L}_{[\a]}=U+\sum_{[\a]\in \Lambda/\sim}I_{[\a]}.
 \end{eqnarray*}
 Finally, Proposition 2.10 gives us $[I_{[\a]},I_{[\b]}]=0$ if $[\a]\neq[\b]$.\hfill $\square$

\begin{definition}
The annihilator of a BiHom-Leibniz superalgebra $\mathfrak{L}$ is the set
\begin{eqnarray*}
Z(\mathfrak{L}):=\{v\in \mathfrak{L}:[v,\mathfrak{L}]+[\mathfrak{L}, v]=0\}.
\end{eqnarray*}
\end{definition}

 \begin{corollary}  If $Z(\mathfrak{L})=0$ and $H=\sum_{\a\in \Lambda}([\mathfrak{L}_{\a\psi^{-1}}, \mathfrak{L}_{-\a\phi^{-1}}])$. Then $\mathfrak{L}$ is the direct sum of the ideals given in Theorem 2.12,
 \begin{eqnarray*}
\mathfrak{L}=\bigoplus_{[\a]\in \Lambda/\sim}I_{[\a]},
\end{eqnarray*}
Furthermore, $[I_{[\a]},I_{[\b]}]=0$ if $[\a]\neq[\b]$.\hfill $\square$.
\end{corollary}
{\bf Proof.} Since $H=\sum_{\a\in \Lambda}([\mathfrak{L}_{\a\psi^{-1}}, \mathfrak{L}_{-\a\phi^{-1}}])$, it follows that  $\mathfrak{L}=\sum_{[\a]\in \Lambda/\sim}I_{[\a]}$. To verify the direct character of the sum, take some $v\in I_{[\a]}\cap(\sum_{[\b]\in\Lambda/\sim,[\b]\neq[\a]}I_{[\b]})$. Since $v\in I_{[\a]}$, the fact $[I_{[\a]},I_{[\b]}]=0$ when  $[\a]\neq[\b]$ gives us
\begin{eqnarray*}
[v, \sum_{[\b]\in\Lambda/\sim,[\b]\neq[\a]}I_{[\b]}]+[\sum_{[\b]\in\Lambda/\sim,[\b]\neq[\a]}I_{[\b]},v]=0.
\end{eqnarray*}
  In a similar way,  $v\in \sum_{[\b]\in\Lambda/\sim,[\b]\neq[\a]}I_{[\b]}$ implies $[v,I_{[\a]}]+[I_{[\a]}, v]=0$. That is $v\in Z(\mathfrak{L})$ and so $v=0$.              \hfill $\square$
\section{The simple components}
\def\theequation{\arabic{section}. \arabic{equation}}
\setcounter{equation} {0}

In this section we focus on the simplicity of split regular BiHom-Leiniz superalgebras $\mathfrak{L}$ by
centering our attention in those of maximal length. From now on $char(\mathbb{K})=0$.

For an ideal  $I$, we have
\begin{eqnarray*}
I=(I\cap H)\oplus (\bigoplus_{\a\in \Lambda} (I\cap \mathfrak{L}_{\a})).
\end{eqnarray*}

\begin{definition}
A split regular BiHom-Leiniz superalgebra $\mathfrak{L}$ is of maximal length if dim$\mathfrak{L}_{\a, \bar{i}}=1$ for any $\a\in \Lambda$ and $\bar{i}\in \mathbb{Z}_2$.
\end{definition}

Observe that if $\mathfrak{L}$ is of  maximal length, then we have
\begin{eqnarray}
I=(I\cap H)\oplus (\bigoplus_{\a\in \Lambda^{I}_{\bar{0}}}  \mathfrak{L}_{\a, \bar{0}}))\oplus  (\bigoplus_{\b\in \Lambda^{I}_{\bar{1}}}  \mathfrak{L}_{\b, \bar{1}})),
\end{eqnarray}
where $\Lambda^{I}_{\bar{i}}=\{\g\in \Lambda: I_{\bar{i}}\cap \mathfrak{L}_{\g, \bar{i}}\neq 0\}, \bar{i}\in \mathbb{Z}_2$.

In particular, case $I=\mathfrak{J}$, we get
\begin{eqnarray}
\mathfrak{J}=(\mathfrak{J}\cap H)\oplus (\bigoplus_{\a\in \Lambda^{\mathfrak{J}}_{\bar{0}}}  \mathfrak{L}_{\a, \bar{0}}))\oplus  (\bigoplus_{\b\in \Lambda^{\mathfrak{J}}_{\bar{1}}}  \mathfrak{L}_{\b, \bar{1}}))
\end{eqnarray}
with $\Lambda^{\mathfrak{J}}_{\bar{i}}=\{\g\in \Lambda: \mathfrak{J}_{\bar{i}}\cap \mathfrak{L}_{\g, \bar{i}}\neq 0\}=\{\g\in \Lambda: 0\neq \mathfrak{L}_{\g, \bar{i}}\subset \mathfrak{J}_{\bar{i}}\}, \bar{i}\in \mathbb{Z}_2$.

From here, we can write
\begin{eqnarray}
\Lambda=\Lambda_{\bar{0}}\cup \Lambda_{\bar{1}}=(\Lambda^{\mathfrak{J}}_{\bar{0}}\cup \Lambda^{\neg\mathfrak{J}}_{\bar{0}})\cup(\Lambda^{\mathfrak{J}}_{\bar{1}}\cup \Lambda^{\neg\mathfrak{J}}_{\bar{1}}),
\end{eqnarray}
where \begin{eqnarray*}
    \Lambda^{\neg\mathfrak{J}}_{\bar{i}}=\{\g\in \Lambda: \mathfrak{L}_{\g, \bar{i}}\neq 0 ~~\mbox{and}~~ \mathfrak{J}_{\bar{i}}\cap \mathfrak{L}_{\g, \bar{i}}= 0\}, \bar{i}\in \mathbb{Z}_2.
      \end{eqnarray*}
      We also denote
      \begin{eqnarray*}
  \Lambda^{\Upsilon}=\Lambda^{\Upsilon}_{\bar{0}}\cup \Lambda^{\Upsilon}_{\bar{1}}, \Upsilon\in \{\mathfrak{J}, \neg\mathfrak{J}\}.
      \end{eqnarray*}
      Therefore, we can write
      \begin{eqnarray}
 \mathfrak{L}=(H_{\bar{0}}\oplus H_{\bar{1}}) \oplus (\bigoplus_{\a\in \Lambda^{\mathfrak{J}}_{\bar{0}}}  \mathfrak{L}_{\a, \bar{0}}))\oplus (\bigoplus_{\b\in \Lambda^{\neg\mathfrak{J}}_{\bar{0}}}  \mathfrak{L}_{\b, \bar{0}}))\oplus (\bigoplus_{\g\in \Lambda^{\mathfrak{J}}_{\bar{1}}}  \mathfrak{L}_{\g, \bar{1}}))\oplus  (\bigoplus_{\delta\in \Lambda^{\neg\mathfrak{J}}_{\bar{1}}}  \mathfrak{L}_{\delta, \bar{1}})).
      \end{eqnarray}

\begin{example}  Consider the 5-dimensional split regular BiHom-Leibniz superalgebra
\begin{eqnarray*}
 \mathfrak{L}=H\oplus \mathfrak{L}_{\a}\oplus \mathfrak{L}_{-\a}\oplus \mathfrak{L}_{\b}\oplus \mathfrak{L}_{-\b}
\end{eqnarray*}
given in Example 2.3, where $ \mathfrak{L}_{\bar{0}}=H\oplus \mathfrak{L}_{\a}\oplus \mathfrak{L}_{-\a}$ and $ \mathfrak{L}_{\bar{1}}=\mathfrak{L}_{\b}\oplus \mathfrak{L}_{-\b}$. This is a split Hom-Leibniz superalgebra  of maximal length such that $\mathfrak{J}=<e_1,e_2>$. Hence $\Lambda^{\neg\mathfrak{J}}_{\bar{0}}={\pm \a}$, $\Lambda^{\mathfrak{J}}_{\bar{1}}={\pm \b}$ and $\Lambda^{\mathfrak{J}}_{\bar{0}}=\Lambda^{\neg \mathfrak{J}}_{\bar{1}}=\emptyset$.
  \end{example}

We are going to refine the concept of connections of nonzero roots in the setup of maximal
length split regular BiHom-Leibniz superalgebras.  We recall that a roots system $\Lambda^{\Upsilon}$ with $\Upsilon\in \{\mathfrak{J}, \neg\mathfrak{J}\}$ of a split regular BiHom-Leiniz superalgebra $\mathfrak{L}$ is called symmetric if it satisfies that $\alpha \in \Lambda^{\Upsilon}$ implies $-\alpha\in \Lambda^{\Upsilon}$. From now on we will suppose that $\Lambda^{\Upsilon}$ with $\Upsilon\in \{\mathfrak{J}, \neg\mathfrak{J}\}$ is symmetric.

\begin{definition}
  Let $\a\in \Lambda^{\Upsilon}_{\bar{i}}$ and  $\b\in \Lambda^{\Upsilon}_{\bar{j}}$   with $\Upsilon\in \{\mathfrak{J}, \neg\mathfrak{J}\}$ and $\bar{i}, \bar{j}\in \mathbb{Z}_2$. We say $\a$ is $\neg\mathfrak{J}$-connected to $\b$, denoted by $\a\sim_{\neg\mathfrak{J}}\b$, if  there exists a family of nonzero roots $\a_1,\a_2,..., \a_n$ such that
\begin{eqnarray*}
\a_k\in \Lambda^{\neg\mathfrak{J}}_{\bar{i_k}}
\end{eqnarray*}
for some $\bar{i_k}\in \mathbb{Z}_2$ and for any $k=2,...,n$; and such that \\
1. $\a_{1}\in \{\a \phi^{-n}\psi^{-r}:~~n, r\in \mathbb{N}\}$.\\
2. $\a_1\phi^{-1}+\a_2\psi^{-1}\in \Lambda^{\Upsilon}_{\bar{i}+\bar{i_2}}$,

 $\a_1\phi^{-2}+\a_2\phi^{-1}\psi^{-1}+\a_3\psi^{-1}\in \Lambda^{\Upsilon}_{\bar{i}+\bar{i_2}+\bar{i_3}}$,

 $\a_1\phi^{-3}+\a_2\phi^{-2}\psi^{-1}+\a_3\phi^{-1}\psi^{-1}+\a_4\psi^{-1}\in \Lambda^{\Upsilon}_{\bar{i}+\bar{i_2}+\bar{i_3}+\bar{i_4}}$,

 $\cdot\cdot\cdot\cdot\cdot$

 $\a_1\phi^{-n+1}+\a_2\phi^{-n}\psi^{-1}+\a_3\phi^{-n+1}\psi^{-1}+\cdot\cdot\cdot+\a_{n-1}\phi^{-1}\psi^{-1}+\a_{n}\psi^{-1}\in \Lambda^{\Upsilon}_{\bar{i}+\bar{i_2}+...+\bar{i_{n}}}$.\\
 3. $\a_1\phi^{-n+1}+\a_2\phi^{-n}\psi^{-1}+\a_3\phi^{-n+1}\psi^{-1}+\cdot\cdot\cdot+\a_{i}\phi^{-n+i}\psi^{-1}+\cdot\cdot\cdot+\a_{n-1}\phi^{-1}\psi^{-1}+\a_{n}\psi^{-1}\in \{\pm \b\phi^{-m}\psi^{-s}:m, s\in \mathbb{N}\}$
 and $\bar{i}+\bar{i_2}+...+\bar{i_{n}}=\bar{j}$.

We also say that $\{\a_1,\a_2,...,\a_k\}$ is a  $\neg\mathfrak{J}$-connected from $\a$ to $\b$.
\end{definition}

Let us introduce the notion of root-multiplicativity in the framework of split
regular BiHom-Leibniz superalgebras of maximal length, in a similar way to the ones for
split regular BiHom-Lie superalgebras in \cite{Zhang2018}.
\begin{definition}
We say that a split regular BiHom-Leibniz superalgebra of maximal length $\mathfrak{L}$ is root-multiplicative if
the below conditions hold.

1. Given $\a\in \Lambda^{\neg\mathfrak{J}}_{\bar{i}}$,  $\b\in \Lambda^{\neg\mathfrak{J}}_{\bar{j}}$ such that $\a\phi^{-1}+\b\psi^{-1}\in \Lambda$ then $[\mathfrak{L}_{\a, \bar{i}}, \mathfrak{L}_{\b, \bar{j}}]\neq 0.$

2.  Given $\a\in \Lambda^{\neg\mathfrak{J}}_{\bar{i}}$,  $\g\in \Lambda^{\mathfrak{J}}_{\bar{j}}$ such that $\a\phi^{-1}+\g\psi^{-1}\in \Lambda^{\mathfrak{J}}$ then $[\mathfrak{L}_{\g, \bar{j}}, \mathfrak{L}_{\a, \bar{i}}]\neq 0.$
\end{definition}
\begin{proposition}
 Suppose $H=\sum_{\a\in \Lambda^{\neg\mathfrak{J}}}([\mathfrak{L}_{\a\psi^{-1}}, \mathfrak{L}_{-\a\phi^{-1}}])$ and $\mathfrak{L}$ is root-multiplicative. If  $\Lambda^{\neg\mathfrak{J}}$ has all of its roots $\neg\mathfrak{J}$-connected,  then any ideal $I$ of $\mathfrak{L}$ such that $I\nsubseteq H\oplus \mathfrak{J}$, then $I=\mathfrak{L}$.
\end{proposition}
{\bf Proof.} By (3.1) and (3.3), we can write
\begin{eqnarray*}
I=(I\cap H)\oplus (\bigoplus_{\a\in \Lambda^{\neg\mathfrak{J},I}_{\bar{0}}}  \mathfrak{L}_{\a, \bar{0}}))\oplus  (\bigoplus_{\b\in \Lambda^{\mathfrak{J}, I}_{\bar{1}}}  \mathfrak{L}_{\b, \bar{1}})),
\end{eqnarray*}
where $\Lambda^{\neg\mathfrak{J},I}=\Lambda^{\neg\mathfrak{J}}\cap \Lambda^{I}$ and $\Lambda^{\mathfrak{J},I}=\Lambda^{\mathfrak{J}}\cap \Lambda^{I}$. Since $I\nsubseteq H\oplus \mathfrak{J}$, there exists $\a_0\in \Lambda^{\neg\mathfrak{J}}_{\bar{i_0}}$ such that
\begin{eqnarray}
0\neq \mathfrak{L}_{\a_0, \bar{i_0}} \subset I, ~~\mbox{for ~some}~ \bar{i_0}\in \mathbb{Z}_2.
\end{eqnarray}
By Lemma 2.6, $\phi(\mathfrak{L}_{\a_0})=\mathfrak{L}_{\a_0\phi^{-1}}$ and $\psi(\mathfrak{L}_{\a_0})=\mathfrak{L}_{\a_0\psi^{-1}}$. (3.5) gives us $\phi(\mathfrak{L}_{\a_0, \bar{i_0}}) \subset \phi(I)=I$ and $\psi(\mathfrak{L}_{\a_0, \bar{i_0}}) \subset \psi(I)=I$. So $\mathfrak{L}_{\a_0\phi^{-1}\psi^{-1}, \bar{i_0}}\subset I$. Similarly we get
\begin{eqnarray}
\mathfrak{L}_{\a_0\phi^{-n}\psi^{-r}, \bar{i_0}}\subset I,  ~~\mbox{for}~ n, r\in \mathbb{N}.
\end{eqnarray}
For any $\b\in \Lambda^{\neg\mathfrak{J}}$, $\b\notin \pm \a_0\phi^{-n+1}\psi^{-1}$, for  $n\in \mathbb{N}$, the fact that $\a_0$ and $\b$ are $\neg\mathfrak{J}$-connected gives us a $\neg\mathfrak{J}$-connection $\{\g_1,\g_2,...,\g_n\}\subset \Lambda^{\neg\mathfrak{J}}$ from $\a_0$ to $\b$ such that

$\g_1=\a_0\in\Lambda^{\neg\mathfrak{J}}_{\bar{i_0}},  \g_k\in \Lambda^{\neg\mathfrak{J}}_{\bar{i_k}}$, for  $k=2,...,n$, \\
 $\g_1\phi^{-1}+\g_2\psi^{-1}\in \Lambda^{\Upsilon}_{\bar{i}+\bar{i_2}}$,~~~ $\g_1\phi^{-n+1}+\g_2\phi^{-n}\psi^{-1}+\g_3\phi^{-n+1}\psi^{-1}+\cdot\cdot\cdot+\g_{n-1}\phi^{-1}\psi^{-1}+\g_{n}\psi^{-1}\in \Lambda^{\Upsilon}_{\bar{i}+\bar{i_2}+...+\bar{i_{n}}}$.\\
 $\g_1\phi^{-n+1}+\g_2\phi^{-n}\psi^{-1}+\g_3\phi^{-n+1}\psi^{-1}+\cdot\cdot\cdot+\g_{i}\phi^{-n+i}\psi^{-1}+\cdot\cdot\cdot+\g_{n-1}\phi^{-1}\psi^{-1}+\g_{n}\psi^{-1}\in \{\pm \b\phi^{-m}\psi^{-s}:m, s\in \mathbb{N}\}$
 and $\bar{i}+\bar{i_2}+...+\bar{i_{n}}=\bar{j}$.

 Taking into account that $\g_1,\g_2$ and $\g_1\phi^{-1}+\g_2\psi^{-1}$.  Since $\g_1=\a_0\in\Lambda^{\neg\mathfrak{J}}_{\bar{i_0}}$ and $\g_2\in \Lambda^{\neg\mathfrak{J}}_{\bar{i_2}}$, the root-multiplicativity
and maximal length of $\mathfrak{L}$ show  $\g_1\phi^{-1}+\g_2\psi^{-1}\in \Lambda_{\bar{i}+\bar{j}}$.
 \begin{eqnarray*}
0\neq [\mathfrak{L}_{\g_1, \bar{i_0}}, \mathfrak{L}_{\g_2, \bar{i_2}}]=\mathfrak{L}_{\g_1\phi^{-1}+\g_2\psi^{-1}, \bar{i_0}+\bar{i_2}}.
 \end{eqnarray*}
 and by (3.6), we have
 \begin{eqnarray*}
0\neq\mathfrak{L}_{\g_1\phi^{-1}+\g_2\psi^{-1}, \bar{i_0}+\bar{i_2}}\subset I.
 \end{eqnarray*}
 We can argue in a similar way from $\g_1\phi^{-1}+\g_2\psi^{-1}, \g_3$ and $\g_1\phi^{-2}+\g_2\phi^{-1}\psi^{-1}+\g_3\psi^{-1}$. Hence
  \begin{eqnarray*}
0\neq [\mathfrak{L}_{\g_1\phi^{-1}+\g_2\psi^{-1}, \bar{i_0}+\bar{i_2}}, \mathfrak{L}_{\g_3, \bar{i_3}}]=\mathfrak{L}_{\g_1\phi^{-2}+\g_2\phi^{-1}\psi^{-1}+\g_3\psi^{-1}, \bar{i_0}+\bar{i_2}+\bar{i_3}}.
 \end{eqnarray*}
 and by above, we have
 \begin{eqnarray*}
0\neq\mathfrak{L}_{\g_1\phi^{-2}+\g_2\phi^{-1}\psi^{-1}+\g_3\psi^{-1}, \bar{i_0}+\bar{i_2}+\bar{i_3}}\subset I.
 \end{eqnarray*}
 Following this process with the $\neg\mathfrak{J}$-connection $\{\g_1,\g_2,...,\g_n\}$ we obtain that
  \begin{eqnarray*}
0\neq\mathfrak{L}_{\g_1\phi^{-n+1}+\g_2\phi^{-n}\psi^{-1}+\g_3\phi^{-n+1}\psi^{-1}+...+\g_n\psi^{-1}, \bar{i_0}+\bar{i_2}+\bar{i_3}+...+\bar{i_n}}\subset I.
 \end{eqnarray*}
 and so that we get  either
 \begin{eqnarray}
\mathfrak{L}_{\b\phi^{-m}\psi^{-s}, \bar{j}}\subset I ~~\mbox{or}~~\mathfrak{L}_{-\b\phi^{-m}\psi^{-s}, \bar{j}}\subset I
 \end{eqnarray}
 for any $\b\in \Lambda^{\neg\mathfrak{J}}, m\in \mathbb{N}$. Moreover, we have
 \begin{eqnarray}
\b\phi^{-m}\psi^{-s}\in \Lambda^{\neg\mathfrak{J}}.
 \end{eqnarray}
Since $H=\sum_{\a\in \Lambda^{\neg\mathfrak{J}}}([\mathfrak{L}_{\a\psi^{-1}}, \mathfrak{L}_{-\a\phi^{-1}}])$, by (3.7) and (3.8), we get
\begin{eqnarray}
H\subset I.
\end{eqnarray}
Now, for any $\Upsilon\in \{\mathfrak{J}, \neg\mathfrak{J}\}$ and $\bar{k}\in \mathbb{Z}_2$, given any $\delta\in \Lambda^{\mathfrak{J}}_{\bar{k}}$, the facts $\delta\neq0, H\subset I$ and the maximal length of $\mathfrak{L}$ show that
\begin{eqnarray*}
\mathfrak{L}_{\delta, \bar{k}}=[\mathfrak{L}_{\delta\psi, \bar{k}}, H_{\bar{0}}]\subset I.
\end{eqnarray*}
The decomposition of $\mathfrak{L}$ in (3.4) finally gives us $H=I$. \hfill $\square$
\begin{definition}
The Lie-annihilator of a split BiHom-Leibniz superalgebra of maximal
length $\mathfrak{L}$ is the set
\begin{eqnarray*}
Z_{Lie}(\mathfrak{L})=\{v \in  \mathfrak{L}: [v, \mathfrak{L}_\a]+[\mathfrak{L}_\a, v]=0, \forall \a\in \Lambda^{\mathfrak{J}}\}.
\end{eqnarray*}
\end{definition}
Observe that $Z(\mathfrak{L})\subset Z_{Lie}(\mathfrak{L})$.
\begin{proposition}
 Suppose $H=\sum_{\a\in \Lambda^{\neg\mathfrak{J}}}([\mathfrak{L}_{\a\psi^{-1}}, \mathfrak{L}_{-\a\phi^{-1}}])$, $Z_{Lie}(\mathfrak{L})=0$ and $\mathfrak{L}$ is root-multiplicative. If  $\Lambda^{\neg\mathfrak{J}}$ has all of its roots $\neg\mathfrak{J}$-connected, then any ideal $I$ of $\mathfrak{L}$ such that $I\subseteq  \mathfrak{J}$, then $I=\mathfrak{J}$.
\end{proposition}
{\bf Proof.} By (3.1), we can write
\begin{eqnarray}
I=(I\cap H)\oplus (\bigoplus_{\b\in \Lambda^{I}_{\bar{0}}}  \mathfrak{L}_{\b, \bar{0}}))\oplus  (\bigoplus_{\g\in \Lambda^{I}_{\bar{1}}}  \mathfrak{L}_{\g, \bar{1}})),
\end{eqnarray}
where
\begin{eqnarray*}
\Lambda^{I}_{\bar{i}}=\{\delta\in \Lambda: I_{\bar{i}}\cap \mathfrak{L}_{\delta, \bar{i}}\neq 0\}=\{\delta\in \Lambda: 0\neq \mathfrak{L}_{\delta, \bar{i}}\subset I_{\bar{i}}\}
\end{eqnarray*}
and with $\Lambda^{I}_{\bar{i}}\subset \Lambda^{\mathfrak{J}}_{\bar{i}}$ for any $\bar{i}\in \mathbb{Z}_2$. Fixed some $\bar{i}\in \mathbb{Z}_2$ and for any $\a\notin \Lambda^{\mathfrak{J}}$ and $\bar{j}\in \mathbb{Z}_2$, we have
\begin{eqnarray*}
[\mathfrak{J}\cap H_{\bar{i}},  \mathfrak{L}_{\a, \bar{j}}]+[ \mathfrak{L}_{\a, \bar{j}}, \mathfrak{J}\cap H_{\bar{i}}]\subset \mathfrak{L}_{\a, \bar{i}+\bar{j}}\subset \mathfrak{J}.
\end{eqnarray*}
Hence, in case $[\mathfrak{J}\cap H_{\bar{i}},  \mathfrak{L}_{\a, \bar{j}}]+[ \mathfrak{L}_{\a, \bar{j}}, \mathfrak{J}\cap H_{\bar{i}}]\neq 0$ we have $\a\in \Lambda^{\mathfrak{J}}$, a contradiction. Hence  $[\mathfrak{J}\cap H_{\bar{i}},  \mathfrak{L}_{\a, \bar{j}}]+[ \mathfrak{L}_{\a, \bar{j}}, \mathfrak{J}\cap H_{\bar{i}}]= 0$ for each $\bar{i}\in \mathbb{Z}_2$, and so
\begin{eqnarray}
\mathfrak{J}\cap H \subset Z_{Lie}(\mathfrak{L}).
\end{eqnarray}
By (3.11), we can also write
\begin{eqnarray}
\mathfrak{J}= (\bigoplus_{\a\in \Lambda^{\mathfrak{J}}_{\bar{0}}}  \mathfrak{L}_{\a, \bar{0}}))\oplus  (\bigoplus_{\delta\in \Lambda^{\mathfrak{J}}_{\bar{1}}}  \mathfrak{L}_{\delta, \bar{1}})).
\end{eqnarray}
Taking into account $I\cap H\subset \mathfrak{J}\cap H=0$, we also write
\begin{eqnarray*}
I=(\bigoplus_{\b\in \Lambda^{I}_{\bar{0}}}  \mathfrak{L}_{\b, \bar{0}}))\oplus  (\bigoplus_{\g\in \Lambda^{I}_{\bar{1}}}  \mathfrak{L}_{\g, \bar{1}})),
\end{eqnarray*}
with $\Lambda^{I}_{\bar{i}}\subset \Lambda^{\mathfrak{J}}_{\bar{i}}$. Hence, we can take some $\b_0\in \Lambda^{I}_{\bar{i}}$ such that
\begin{eqnarray*}
0\neq \mathfrak{L}_{\b_0, \bar{i}}\subset I.
\end{eqnarray*}
Now, we can argue with the root-multiplicativity and the maximal length of $\mathfrak{L}$ as in Proposition 3.5 to
conclude that given any $\b \in \Lambda^{\mathfrak{J}}_{\bar{j}}, \bar{j}\in \mathbb{Z}_2$, there exists a $\neg\mathfrak{J}$-connection
$\{\delta_1,\delta_2, ..., \delta_r\}$ from $\b_0$ to $\b$ such that
\begin{eqnarray*}
0\neq [[...[\mathfrak{L}_{\b_0, \bar{i}}, \mathfrak{L}_{\delta_2, \bar{i_2}}],...], \mathfrak{L}_{\delta_r, \bar{i_r}}]\in \mathfrak{L}_{\b\phi^{-m}\psi^{-s}, \bar{j}}, ~\mbox{for}~m, s\in \mathbb{N}
\end{eqnarray*}
and so
\begin{eqnarray}
\mathfrak{L}_{\epsilon\b\phi^{-m}\psi^{-s}, \bar{j}}\subset I, ~~\mbox{for~some}~\epsilon\in \pm1, m, s\in \mathbb{N}.
\end{eqnarray}
Note that $\b\in \Lambda^{\mathfrak{J}}_{\bar{j}}$ indicates $\mathfrak{L}_{\b, \bar{j}}\in \mathfrak{J}$.   By Lemma 2.6, $\phi(\mathfrak{L}_{\b})=\mathfrak{L}_{\b\phi^{-1}}$ and $\psi(\mathfrak{L}_{\b})=\mathfrak{L}_{\b\psi^{-1}}$. Since $\mathfrak{L}$ is of maximal length, we have  $\phi(\mathfrak{L}_{\b, \bar{j}}) \subset \phi(\mathfrak{J})=\mathfrak{J}$ and $\psi(\mathfrak{L}_{\b, \bar{j}}) \subset \psi(\mathfrak{J})=\mathfrak{J}$. So $\mathfrak{L}_{\b\phi^{-1}\psi^{-1}, \bar{j}}\subset I$. Similarly we get
\begin{eqnarray}
\mathfrak{L}_{\b\phi^{-m}\psi^{-s}, \bar{j}}\in \mathfrak{J}, ~\mbox{for}~m, s\in \mathbb{N}.
\end{eqnarray}
Hence we can argue as above with the
root-multiplicativity and maximal length of $\mathfrak{L}$ from $\b$ instead of $\b_0$, to get that in case $\epsilon \b_0\phi^{-m}\psi^{-s}\in \Lambda^{\mathfrak{J}}_{\bar{k}}$
for some $\epsilon\in  \pm1$ and $\bar{k}\in \mathbb{Z}_2$, then $0\neq \mathfrak{L}_{\epsilon\b_0\phi^{-m}\psi^{-s}, \bar{k}}\in I$.

The decomposition of $\mathfrak{J}$ in (3.12) finally gives us $I=\mathfrak{J}$.
\hfill $\square$

\begin{theorem} Let $\mathfrak{L}$ be a  split regular BiHom-Leiniz superalgebra  of maximal length,  $Z_{Lie}(\mathfrak{L})=0$. and $\mathfrak{L}$ is root  multiplicative. Then $\mathfrak{L}$ is simple if and only if it is prime and $\Lambda^{\mathfrak{J}}, \Lambda^{\neg\mathfrak{J}}$ has all of their roots $\neg\mathfrak{J}$-connected.
\end{theorem}
{\bf Proof.} Consider  any $\a\in \Lambda^{\neg\mathfrak{J}}_{\bar{i}}$  with $i\in \mathbb{Z}_2$ and the subspace $\mathfrak{L}_{\a, \bar{i}}$. Let us denote
by $I(\mathfrak{L}_{\a, \bar{i}})$ the ideal of $\mathfrak{L}$ generated by $\mathfrak{L}_{\a, \bar{i}}$.  By simplicity $I(\mathfrak{L}_{\a, \bar{i}})=\mathfrak{L}$. We observe that the fact $\mathfrak{J}$ is an ideal of $\mathfrak{L}$ and
we assert that $\mathfrak{L}_{\a, \bar{i}}$ is contained in the linear span of the set
\begin{eqnarray*}
&&\{[[...[v_\a, v_{\b_1}], ...,],v_{\b_n}]; [v_{\b_n},[...[v_{\b_1}, v_\a],...]];\\
&&[[...[v_{\b_1}, v_\a],...], v_{\b_n};[v_{\b_n},[...[v_\a, v_{\b_1}],...]] ~~\mbox{with}~~0\neq v_{\a}\in \mathfrak{L}_{\a, \bar{i}},\\
&&0\neq v_{\b_i}\in \mathfrak{L}_{\b_i, \bar{j_i}}, \b_i\in \Lambda,  \bar{j_i}\in \mathbb{Z}_2, n\in \mathbb{N}\}
\end{eqnarray*}
From here, given any $\a'\in \Lambda^{\neg\mathfrak{J}}_{\bar{j}}$, the above observation
gives us that we can write $\a'=\a\phi^{-n}+\b_1\phi^{-n+1}\psi^{-1}+...+\b_n\psi^{-1}$ with any $\b_i\in \Lambda^{\neg\mathfrak{J}}_{\bar{j_i}}$
and being the partial sums $\a\phi^{-n}+\b_1\phi^{-n+1}\psi^{-1}+...+\b_k\phi^{-n+k}\psi^{-1}\in \Lambda^{\neg\mathfrak{J}}_{\bar{i}+\bar{j_1}+...,+\bar{j_k}}$. From here, we have that $\{\a, \b_1,\b_2,...,\b_n\}$
is a $\neg\mathfrak{J}$-connection from $\a$ to $\a'$  and we can assert that $ \Lambda^{\neg\mathfrak{J}}$ has all of its elements $\neg\mathfrak{J}$-connected.

If $\Lambda^{\mathfrak{J}}\neq \emptyset$ and we take some $\b\in \Lambda^{\mathfrak{J}}_{\bar{i}}$ with $i\in \mathbb{Z}_2$,  a similar above argument gives us $ \Lambda^{\neg\mathfrak{J}}$ has all of its elements $\neg\mathfrak{J}$-connected.

The converse is a consequence of Proposition 3.5 and Proposition 3.7.   \hfill $\square$

 \begin{center}
 {\bf ACKNOWLEDGEMENT}
 \end{center}

 The paper is supported by  the NSF of China (No. 11761017),  the Youth Project for Natural Science Foundation of Guizhou provincial department of education (No. KY[2018]155).

\end{document}